\documentclass{article}
\usepackage{amsthm}
\usepackage{amsmath}
\usepackage{amssymb}
\usepackage{rotating}

\newcommand{\rot}[1]{\begin{turn}{270}#1\end{turn}}

\usepackage[all,cmtip]{xy}
\usepackage[all]{xy}
\newtheorem*{theorem}{Theorem}
\newtheorem*{lemma}{Lemma}

\newtheorem{definition}{Definition}

\def\O{{\mathcal O}}
\def\Ac{{\mathcal A}}
\def\Bc{{\mathcal B}}
\def\Mc{{\mathcal M}}

\def\Oc{{\mathcal O}}

\def\bz{{\bf z}}

\def\Hom{{\mathcal H om}}
\def\banica{{B\u anic\u a}}

\def\Proof{\noindent\hskip4pt{\it Proof}.\ }
\def\qed{\hfill$\Box$\vskip10pt}

\begin{document}

\title{Cuspidal Multiple Structures on Smooth Algebraic Varieties as Support}
\author{Nicolae Manolache}
\date{}

\maketitle

\begin{flushright}
 This paper is dedicated to \c Serban Basarab on his 70th Anniversary
\end{flushright}

\section{Introduction}

The aim of this paper is to describe two new classes of locally complete intersection
(lci for short) nilpotent structures on a smooth algebraic variety as support, which we
call {\em cuspidal of types $C_{2,n}$, $C_{3,n}$}. We recall the  known classes: 1)
``primitive or ``quasiprimitive'' structures constructed in \cite{BF1}, \cite{BF2}, and
studied by several authors (cf. \cite{Bo},
\cite{Dr} ) (given locally -- respectively in the general point -- by  ideals of the form
$(x^n, z_1, \ldots z_r)$ ) and 2) ``next`` locally monomial case, constructed in
\cite{M4},
where a class of multiple structures which contains the lci multiple
structures defined locally by an ideal of the shape $(x^n,y^2 , z_1 \ldots z_r)$ is
studied.

The two classes of lci structures on a smooth algebraic variety as support studied here  
are characterized by ideals which locally have the form $(y^2+x^n , xy ,
z_1, \ldots , z_r)$ respectively  $(y^3+x^n , xy , z_1, \ldots , z_r)$, in convenient
local parameters $x$, $y$, $z_1, \ldots z_r$.

\section{Preliminaries}

Let $X$ be a smooth connected algebraic variety over an algebraically closed field $k$ and
a locally Cohen-Macaulay scheme $Y$ such that $Y_{red}$ is $X$. In this case $Y$ is called
a \emph{multiple structure on $X$} and all local rings of $Y$ have the same
multiplicity (cf. \cite{M1}), which is called \emph{the multiplicity of $Y$}. Let $Y$ be
embedded in a smooth variety $P$.
To  $Y$ one associates canonically three filtrations. 
Let $I$ be the (sheaf) ideal of $X$ in $P$ and $J$ be the ideal of $Y$ in $P$. 
Let $m$ be the positive integer such that $I^m \not\subset J$,
$I^{m+1} \subset J$. The three filtrations are:\\
1. Let $I^{(\ell)}$ be the ideal obtained throwing away the embedded components of
$I^\ell+J$ and let 
$Z_\ell$ be the corresponding scheme. This  gives the 
\emph{\banica -Forster filtration} (cf. \cite{BF2}):
\begin{displaymath}
\begin{array}{ccccccccccccc}
 \O_Y=I^{(0)} & \supset & I= I^{(1)} & \supset &  I^{(2)} & \supset &  \ldots &  \supset
&   I^{(m)} & 
\supset & I^{(m+1)}=0\medskip \\
&  & X=Z_1 & \subset & Z_2 & \subset &  \ldots &  \subset & Z_m & \subset & Z_{m+1} =  Y
\end{array}
\end{displaymath}
$Z_\ell$ are not, in general, Cohen-Macaulay. But this is true if ${\rm dim}(X)=1$.
The graded associated object $\Bc (Y)=\bigoplus_{\ell =0}^m I^{(\ell)} / I^{(\ell +1)} $
is naturally a 
graded $\O_X$-algebra. If the schemes $Z_\ell$ are Cohen-Macaulay, the graded components
of $\Bc (Y)$ are locally 
free sheaves on $X$.\\
2. Let $X_\ell $ be defined by $I_\ell =J:I^{m+1-\ell}$. Again, if ${\rm dim}(X)=1$,
$X_\ell $ are Cohen-Macaulay.
This is also true if $Y$ is lci (i.e. locally complete intersection) of multiplicity
at most $6$ (cf. \cite{M2}).
In general this is not always the case. When $X_\ell$ are Cohen-Macaulay, the quotients
$I_\ell /I_{\ell +1}$ are 
locally free sheaves on $X$. This filtration was considered in \cite{M1}.\\
3. Let $Y_\ell$ be the scheme given by $J_\ell=J:I_{m+1-\ell}=J:(J:I^\ell)$. When
$X_\ell$ is Cohen-Macaulay, $Y_\ell$
has the same property. The graded object $\Ac (Y)=\bigoplus_{\ell =0}^m J_\ell
/ J_{\ell +1}$ is a graded
$\O _X$-algebra and $\Mc (Y) =\bigoplus_{\ell =0}^m I_\ell / I_{\ell+1}$ is a graded
$\Ac(Y)$-module.
This filtration was considered in \cite{M2}.\\
The system of the graded components ($\Ac_0(Y),\ldots  \Ac_m(Y);\Mc_0(Y),\ldots
\Mc_m(Y))$ is called 
\emph{the type of $Y$}. $Y$ is called \emph{of free type} when all the graded pieces
are locally free.
As already remarked, in dimension $1$, or if $Y$ is lci of multiplicity up to $6$, this
is  the case.

Recall some properties:\\
1) In general the above {\bf filtrations are different}.
Take for instance $X=Spec(k)$, $Y=\hbox{Spec}(k[x,y]/(x^3,xy,y^4))$,
$P=\hbox{Spec}k[x,y]$\\
2) $Z_\ell \subset Y_\ell \subset X_\ell$ \\
2') there are {\bf canonical morphisms}: $\Bc(Y)\to \Ac(Y)\to \Mc (Y)$\\
3) The multiplications
\begin{displaymath}
\begin{array}{lll}
\Ac_{\ell_1}\otimes \Ac_{\ell_2} & \to & \Ac_{\ell_1+\ell_2}\medskip\\
\Ac_{\ell_1}\otimes \Mc_{\ell_2} & \to & \Mc_{\ell_1+\ell_2}
\end{array}
\end{displaymath}
{\bf are never the zero} morphisms for $\ell_1, \ell_1 \ge 0$, $\ell_1+\ell_2 \le m$ (cf.
\cite{M2}.\\
4) There are canonical {\bf edge morphisms} $ \Mc _{m-1} \to \Ac _1$\\
5) One has the {\bf exact sequences}:
\begin{displaymath}
\begin{array}{c}
0\to \Mc_\ell(Y)\to \O_{X_{\ell+1}} \to \O_{X_\ell} \to 0\medskip\\
0\to \Ac_\ell(Y)\to \O_{Y_{\ell+1}} \to \O_{Y_\ell} \to 0
\end{array}
\end{displaymath}
6) If $Y$ is Gorenstein of free type, then $X_\ell $ and $Y_{m+1-\ell}$ are {\bf locally
algebraically linked}
(cf. \cite{M1}. In particular one has the exact sequences:
\begin{displaymath}
\begin{array}{c}
0\to \omega_{X_{m+1-\ell}}\otimes \omega_Y^{-1}\to \O_Y \to \O_{Y_\ell} \to 0\medskip\\
0\to \omega_{Y_{m+1-\ell}}\otimes \omega_Y^{-1}\to \O_Y \to \O_{X_\ell} \to 0
\end{array}
\end{displaymath}
7) {\bf Duality}. Let $Y$ be a free type Cohen-Macaulay multiple structure on a smooth
support $X$.\\
Then $Y$ is Gorenstein if and only if the following conditions are fulfilled:\\
a) $\Ac_m$ and $\Mc_m$ are line bundles\\
(b) $\Ac_m=\Mc_m$  \\
(c) The canonical morphisms:
\begin{displaymath}
\Ac_\ell \to \Hom _{\O_X}(\Mc_{m-\ell},\Mc_m)\cong \Mc_{m-\ell}^\vee\otimes \Mc_\ell
\end{displaymath}
are isomorphisms (cf \cite{M3})\\
7') In particular: if $Y$ is Gorenstein of free type, then (cf also \cite{M2}):
\begin{displaymath}
\begin{array}{l}
(a)\  \rm{rank}\ \Ac_\ell(Y)=\rm{rank}\ \Mc_{m-\ell}(Y)\medskip\\
(b)\  \Ac_\ell(Y)=\Mc _\ell (Y) \rm{\ iff\ } \rm{rank}\ \Ac_\ell(Y)=\rm{rank}\
 \Ac_{m-\ell}(Y)
\end{array}
\end{displaymath}
In this paper all the schemes are algebraic schemes over a fixed algebraically closed
field $k$, 
of characteristic $0$.

\section{Cuspidal Multiple Structures}

\begin{definition}
Let $X$ be  a smooth variety embedded in a smooth one $P$.
Suppose $\hbox{codim}_PX \ge 2$. We say that a nilpotent 
 scheme structure $Y \subset P$ on $X$ 
is a cuspidal nilpotent structure (of type $C_{m,n}$) if, in any point $p \in X$, there
are local parameters such that the completed local rings have the following shape:
\begin{align*}
\widehat{\mathcal O}_{p,X} & \cong  k[[u_1,\ldots ,u_d]], \\
\widehat{\mathcal O}_{p,P} & \cong  k[[u_1, \ldots , u_d, x,y,z_1, \ldots , z_r ]] \\
\widehat{\mathcal O}_{p,Y} & \cong  k[[u_1,\ldots , u_d, x,y,z_1, \ldots ,
z_r]]/(y^m+x^n,xy,
z_1,\ldots ,z_r) \ , m\le n
\end{align*}
\end{definition}

\noindent In the following all local shapes of various ideals are considered in
$\widehat{\mathcal O}_{p,X}$.

\subsection{ $C_{2,n}$}

We assume $n \ge 3$. The case $n=2$ is treated in \cite{M1}, \cite{M2}. In the following
one denotes ${\bf z}=(z_1, \ldots ,z_r)$.
If $Y$ is a cuspidal nilpotent structure of type $C_{2,n}$ on $X \subset P$,
then the canonical filtrations, look locally:
\begin{center}$
\begin{array}{llll}
J : I^0 & = J =I_{n+1}  & J:(J:I^0) & =\Oc \\
J:I^1 &=(x^n,xy,y^2,{\bz} )=I_n & J:(J:I^1) & =(x,y,\bz)=J=J_1\\
J:I^2 & =(x^{n-1},y,{\bz} )=I_{n-1} & J:(J:I^2) & =(x^2,xy,y^2, \bz)=J_2\\
J:I^3 & =(x^{n-2},y,{\bz} )=I_{n-2} & J:(J:I^3) & =(x^3,xy,y^2, \bz)=J_3\\
\hfill \vdots \hfill & & \hfill \vdots \hfill \\
J:I^{n-1} & =(x^2,y,{\bz} )=I_2 & J:(J:I^{n-1}) & =(x^{n-1},xy,y^2, \bz)=J_{n-1}\\
J:I^n& =(x,y,{\bz} ) =I=I_1 & J:(J:I^n) & =(x^n,xy,y^2, \bz)=J_n\\
J:I^{n+1} & =\Oc=I_0 & J:(J:I^{n+1}) & =J=J_{n+1}
\end{array}
$
\end{center}
Let $X_\ell$ be the scheme defined by $I_\ell$ and  $Y_\ell$ the scheme defined
by $J _\ell$.
Then $I_1/I_2 = L$ is a line bundle on $X$ and the scheme $X_2$ defined by $I_2$ is a
double structure on $X$. Observe that $L\otimes L =I/I_2 \otimes I/I_2
\xrightarrow{multiplication} I^2/II_2$ is surjective, so an isomorphism  $L^2 \cong
I^2/II_2$. As the morphism $I^2/II_2 \to I_2/I_3 $ is a surjection between locally free
rank 1 sheaves on $X$, it follows $I_2/I_3 \cong L^2$. In a similar way one shows $I_\ell
/I_{\ell +1} \cong I^\ell/I^{\ell -1}I_2 \cong L^\ell$ for $\ell =3,\ldots ,n-2$.
$E':=I_{n-1}/I_n$ is obviously a rank 2 vector bundle on $X$ and $I_n/I_{n+1} =L^n$. So:
$$
\Mc_Y =\Oc _X \oplus L\oplus \ldots \oplus L^{n-2}\oplus E'\oplus L^n \ \ \ .
$$
Analogously:
$$
\Ac_Y=\Oc _X\oplus E\oplus L^2\oplus \ldots \oplus L^n  \ \ \ ,
$$
where $E=J/J _2$, $L^\ell = J_\ell /J _{l+1}$, $\ell =2,\ldots ,n$.

In almost all degrees, the canonical morphism $\Ac _Y \to \Mc _Y $ is an isomorphism. The
exceptions are the surjection $(\Ac _Y)_1 = E \to L=(\Mc _Y)_1$ and the injection $(\Ac
_Y)_{n-1} = L^{n-1} \to E'=(\Mc _Y)_{n-1}$.
We denote by $K$ the kernel of $E \to L $, so $K=I_2/J _2$. The duality $\Mc _\ell
\cong \Ac _{n-\ell}^\vee \otimes L^n$ gives $E' \cong E^\vee \otimes L^n$. The canonical
morphism $\Mc _{n-1} \to \Ac _1 $ completes to the exact sequence:
\small{
\begin{center}
$
\begin{array}{ccccccccccc}
0 & \hspace{-1em}\to &\hspace{-1em} L^{n-1} &\to \hspace{-1em}& E' & \hspace{-1em}\to
\hspace{-1em}& E & \hspace{-1em}\to \hspace{-1em}& L &\hspace{-1em} \to & 0
\\
 &   & \parallel &   & \parallel &   & \parallel &   & \parallel &  &  \\
0 & \hspace{-1em}\to &\hspace{-1em} \frac{(x^{n-1}, xy, y^2, \bz)}{(x^n, xy, y^2,\bz)}
&\hspace{-1em}\to \hspace{-1em}& \frac{(x^{n-1}, y,
 \bz)}{(x^n, xy, y^2,\bz)}  & \hspace{-1em}\to \hspace{-1em}& \frac{(x, y,  \bz)}{(x^2,
xy, y^2,\bz)} &\hspace{-1em} \to \hspace{-1em}&
\frac{(x, y,  \bz)}{(x^2, xy, y^2,\bz)} & \hspace{-1em}\to & 0

\end{array}
$
\end{center}
}
This decomposes into two exact sequences:
\begin{center}
$
\begin{array}{ccccccccc}
0 & \to & L^{n-1} & \to & E^\vee \otimes L^n  & \to & K &\to & 0 \\
0 & \to & K       & \to & E                   & \to & L & \to & 0
\end{array}
$
\end{center}
and so $K^2 \cong L^n$ .\medskip

\noindent A necessary condition, less obvious, is given  next:

\begin{lemma}
If a cuspidal structure $Y$ of type $C_{2,n}$ does exists on $X \subset P$, with $E$, $L$,
$K$ as above, then one should have an exact sequence:
$$
0 \to L \to E \to K \to 0
$$
\end{lemma}
\Proof As $K \cong I_2/J_2$, the multiplication gives: $K \otimes K \cong I_2/J_2 \otimes
I_2/J_2 \to I_2^2/I_2J_2$ which is a surjective morphism of locally free rank one sheaves
on $X$, hence an isomorphism.

\noindent As the morphism $ I_2^2/I_2J_2 \to
II_2/(I^2 \cap J_3)$ is also a surjective morphism of locally free rank $1$ sheaves on
$X$, it is also an isomorphism. So: $K^2 \cong II_2/(I^2 \cap J_3)$.

\noindent Similarly one shows: $E \otimes K \cong I/J_2 \otimes I_2/J_2 \cong II_2/IJ_2.$

\noindent As $IJ_2 \subset I^2 \cap J_3$ we get the exact sequence:
$$
0 \to \frac{I^2 \cap J_3}{IJ_2} \to \frac{II_2}{IJ_2} \to \frac{II_2}{I^2 \cap J_3} \to 0
$$
It is easy to show that the first nonzero term of this sequence is isomorphic to
$L \otimes K$. Indeed, the multiplication gives : $L \otimes K \cong I/I_2 \otimes
I_2/J_2 \to II_2 /(IJ_2 +I_2^2)$, which must be an isomorphim, being a surjection of
locally free rank $1$ sheaves on $X$. Moreover: $I^2 \cap J_3 \subset II_2$, so that one
has a morphism:
$$
\frac{I^2 \cap J_3}{IJ_2} \to \frac{II_2}{IJ_2 + I_2^2} \ .
$$
which is again a surjection between locally free rank $1$ sheaves on $X$, and so an
isomorphism. So far we get an exact sequence:
$$
0 \to L \otimes K \to E \otimes K \to K^2 \to 0 \ ,
$$
where from the exact sequence of the lemma.
\qed

\noindent Conversely:
\begin{theorem}
Let $X \subset P$ be two smooth varieties, $\hbox{codim}_PX \ge2$. Let $I$ be
the sheaf-ideal of $X$ in $P$. All cuspidal nilpotent structures of type $C_{2,n}$ on $X$,
embedded in $P$ can be obtained in the following way:\\
{\bf Step 1}. Give two line bundles $L$, $K$ on $X$, satisfying $L^n \cong K^2$ and
an extension $0 \to K \stackrel{\iota}{\longrightarrow} E \stackrel{\pi}{\longrightarrow}
L \to 0$  , such that there exists also an extension $0 \to L \to E \to K \to 0$. \\
{\bf Step 2}. Give two surjections $p_2:I/I^2 \to E$, $q_2: I/I^2 \to L$ such that the
diagram :
\begin{center}
$
\begin{array}{ccl}
I/I^2 & \stackrel{p_2}{\longrightarrow} & E \\
\parallel & & \downarrow \pi \\
I/I^2 & \stackrel{q_2}{\longrightarrow} & L
\end{array}
$
\end{center}
is commutative. Take $J_2=\hbox{ker} (I \to I/I^2 \to E)$, $I_2=\hbox{ker} (I \to I/I^2
\to L)$. Then $K \cong I_2/J_2$, $K \otimes E \cong II_2/IJ_2$, and, in convenient
''local coordinates'', $J_2 = (x^2, xy, y^2, \bz)$, $I_2 = (x^2, y, \bz)$. \\
{\bf Step 3}. Give a retract $q_3$ of the canonical inclusion $L^2 \cong I^2/II_2
\hookrightarrow
I_2/II_2$, give $p_3:J_2/IJ_2 \to L^2$ a surjection which makes commutative the diagram

$$
\xymatrix{
J_2/IJ_2  \ar[r]^{p_3} \ar[d] &  L^2 \ar[r] \ar@{=}[d] & 0 \\
I_2/II_2  \ar[r]^{q_3} & L^2 \ar[r] \ar[d]^{\rot{$\cong$}} & 0 \\
& I/II_2 \ar@{_{(}->}[lu] & \\
}
$$
and take $J_3=\hbox{ker}(J_2\to I/IJ_2 \to L^2)$, $I_3=\hbox{ker}(I_2\to I/II_2 \to L^2)$.
Then, in convenient local coordinates, $J_3=(x^3,xy,y^2,\bz)$, $I_3=(x^3,y,\bz)$, and

$K \cong I_3/J_3$, $E\otimes K \cong \frac{II_3}{IJ_3}$.

\hfil \vdots
\medskip

\noindent {\bf Step $\ell$, $\ell =4, \ldots , n-1$}. Give a retract $q_\ell$ of the
canonical inclusion $L^{\ell -1} \cong I^{\ell -1}/I^{\ell -2}I_2 \hookrightarrow I_{\ell
-1}/II_{\ell -1}$,
give $p_\ell : J_{\ell -1}/IJ_{\ell -1} \to L^{\ell -1}$ a surjection which makes
commutative the diagram:
$$
\xymatrix{
J_{\ell -1}/IJ_{\ell -1}  \ar[r]^{p_\ell} \ar[d] &  L^{\ell -1} \ar[r] \ar@{=}[d] & 0 \\
I_{\ell -1}/II_{\ell -1}  \ar[r]^{q_\ell } & L^{\ell -1} \ar[r] \ar[d]^{\rot{$\cong$}} & 0
\\
& I^\ell /I^{\ell -1}I_2 \ar@{_{(}->}[lu] & \\
}
$$\medskip
\noindent
and take $J_\ell =\hbox{ker}(J_{\ell -1} \to J_{\ell -1}/IJ_{\ell -1} \to L^{\ell -1})$,
$I_{\ell}=\hbox{ker}(I_{\ell -1}\to I_{\ell -1}/II_{\ell -1 } \to L^{\ell -1})$.
Then, in convenient ''local coordinates``, $J_\ell =(x^\ell ,xy,y^2,\bz)$,
$I_\ell =(x^\ell ,y,\bz)$, $K\cong I_\ell /J_\ell $
$E\otimes K \cong \frac{II_\ell }{IJ_\ell }
\cong \frac{(x^{\ell +1} ,xy, y^2 , x\bz ,y\bz , \bz ^2)}{(x^{\ell +1} , x^2y, xy^2 ,y^3,
x\bz , y\bz , \bz ^2)}$.

\hfil \vdots
\medskip

\noindent {\bf Step n}. Suppose first  $n \ge 4$. Then give a retract $p_n$ of the
canonical inclusion $L^{n-1} \cong I^n/I^{n-2}I_2 \hookrightarrow J_{n-1}/IJ_{n-1} $. Take
$J_n = I_n =\hbox{ker}(J_{n-1} \to J_{n-1}/IJ_{n-1} \to L^{n-1})$. Then $II_{n-1}\subset
J_n \subset I_{n-1}$, and $J_n/II_{n-1} \to I_{n-1}/II_{n-1}$ is injective. Denote by $E'$
the cokernel of this morphism. Then $E'$ is a vector bundle of rank 2.\\
If $n=3$,  $L^{n-1}=L^2$ is no longer a subbundle of $J_{n-1}/IJ_{n-1}=J_2/IJ_2$. The
morphism $p_3$ is chosen such that the following diagram to be commutative:
$$
\xymatrix{
J_2/IJ_2  \ar[r]^{p_3} \ar[d] &  L^2 \ar[r]  \ar@{^{(}->}[ld]  & 0 \\
I_2/II_2   &  & \\
}
$$\medskip
and then proceed as in the case $n \ge 4$.
In local convenient coordinates $J_n=I_n=(x^n,xy,y^2,\bz)$

\noindent {\bf Step n+1}. The morphism
$$K \otimes K \cong (I_2/J_2) \otimes
(I_{n-1}/J_{n-1}) \xrightarrow{multiplication} I_2I_{n-1}/(I_2J_{n-1}+I_{n-1}J_2)
$$
is a surjection of locally free rank $1$ sheaves, hence an isomorphism.

\noindent Take $p_{n+1} : I_n/II_n \to L^n \cong K^2$ a retract of the
canonical inclusion $L^n \cong I^n/I^{n-1}I_2 \hookrightarrow I_n/II_n$ and of the
canonical inclusion $ K^2 \cong I_2I_{n-1}/(I_2J_{n-1}+I_{n-1}J_2)
\hookrightarrow I_n/II_n$. Then, locally, in convenient coordinates, $J_{n+1} \cong (y^2
+x^n , xy, \bz )$, so that $J_{n+1}$ defines a $C_{2,n}$ cuspidal multiple structure $Y
\subset P$ on  $X$.

\end{theorem}

\Proof All we have to do is to carefully verify, mainly by computation, the assertions
made in the theorem.

\qed

\subsection{$C_{3,n}$}
We assume $n \ge 4$. The case $n=3$ is treated in \cite{M2}.

\noindent If $Y$ is a cuspidal nilpotent structure of type $C_{3,n}$ on $X \subset P$,
then the canonical filtrations look locally:
\begin{center}$
\begin{array}{llll}
J : I^0 &= J=I_{n+1}  & J:(J:I^0) &=\Oc \\
J:I^1 &= (x^n,xy,y^3,{\bz} )=I_n & J:(J:I^1) &=(x,y,\bz)=J=J_1\\
J:I^2 &= (x^{n-1},xy, y^2,{\bz} )=I_{n-1} & J:(J:I^2)&=(x^2,xy,y^2, \bz)=J_2\\
J:I^3 &= (x^{n-2},y,{\bz} )=I_{n-2} & J:(J:I^3)&=(x^3,xy,y^3, \bz)=J_3\\
\hfill \vdots \hfill & & \hfill \vdots \hfill \\
J:I^{n-1} &= (x^2,y,{\bz} )=I_2 & J:(J:I^{n-1})&=(x^{n-1},xy,y^3, \bz)=J_{n-1}\\
J:I^n &= (x,y,{\bz} )=I=I_1 & J:(J:I^n)&=(x^n,xy,y^3, \bz)=J_n\\
J:I^{n+1} &= \Oc=I_0 & J:(J:I^{n+1})&=J=J_{n+1}
\end{array}
$
\end{center}

\noindent Let $X_\ell $ be the scheme defined by $I_\ell$ and $Y_\ell $ the scheme defined
by
$J_\ell$. 

\noindent One shows that the associated graded objects have the shape:
$$
\Mc_Y =\Oc _X \oplus L\oplus L^2 \oplus \ldots \oplus L^{n-3}\oplus F'\oplus E' \oplus L^n
\ \ \ .
$$
Analogously:
$$
\Ac_Y=\Oc _X\oplus E\oplus F \oplus L^3\oplus \ldots \oplus L^n  \ \ \ ,
$$
where  $L^\ell=I_\ell /I_{\ell +1} \cong I^\ell /I^{\ell -1}I_2$, $\ell
=1,\ldots ,n-3 $,  $F'=I_{n-2}/I_{n-1}$, $E' =I_{n-1}/I_n$, $E=I/J_2$,
$F=J_2/J_3$. $L^\ell = J_\ell /J _{l+1}$, $\ell =3,\ldots ,n$.

\noindent The canonical morphism $\Ac _Y \to \Mc _Y$ is an isomorphism in all degrees
except degrees $1$, $2$, $n-2$, $n-1$. We analyze these situations:

\noindent {\em Degree $1$}: The morphism $E \to L$ is an epimorphism, and let $K$ be the
kernel, i.e.
$K:=I_2/J_2$.

\noindent {\em Degree $2$}: One proves easily $S^2E \cong I^2/IJ_2$, $E \otimes K \cong
II_2/IJ_2$, $K^2 \cong I_2^2/I_2J_2$. Also, the canonical morphism $I_2^2/I_2J_2 \to
II_2/I^2 \cap J_3$ is a surjection between two locally free rank $1$ sheaves, i.e. an
isomorphism.
Moreover, the surjectivity of the canonical morphism $II_2/I^2 \cap J_3 \to I_3 \cap
J_2/J_3$ gives a new expression of $K^2$. Summing up: $K^2 \cong I_2^2/I_2J_2 \cong
II_2/I^2 \cap J_3 \cong I_3 \cap J_2/J_3$. The
multiplication in $\Ac_Y$ gives a morphism: $S^2E \to F$, which completes to an exact
sequence:
$$
0 \to I^2 \cap J_3 /IJ_2 \to I^2/IJ_2 \to J_2/J_3 \to 0
$$
Observe that $L \otimes K \cong II_2/(I_2^2 + IJ_2)$, and, as $I^2\cap J_3 \subset II_2$,
one has a canonical morphism $I^2 \cap J_3 /IJ_2 \to II_2/(IJ_2+I_2^2)$, which, as a
surjective morphism of locally free rank $1$ sheaves, is an isomorphism. The above exact
sequence becomes:
$$
0 \to L\otimes K \to S^2E \to F \to 0 \ . 
$$
The canonical surjection $E\to L$ gives the surjection $S^2E \to L^2$, which
compleproducestes to the exact sequence:
$$
0 \to II_2/IJ_2 \to I^2/IJ_2 \to I^2/II_2 \to 0
$$
i.e.
$$
0 \to E\otimes K \to S^2E \to L^2 \to 0 \ .
$$

\noindent One has the exact sequence:
$$
0 \to I_3 \cap J_2/J_3 \to J_2/J_3 \to I_2/I_3 \to 0 \ ,
$$
which translates to:
$$
0 \to K^2 \to F \to L^2 \to 0 .
$$

\noindent This fits in the commutative diagram:
$$
\xymatrix{
& & S^2E \ar[d] \ar[rd] & & \\
0 \ar[r] & K^2 \ar[r] \ar@{^{(}->}[ru] &  F  \ar[r]  &
L^2 \ar[r]  & 0 \\
}
$$

\noindent {\em Degree $n-2$}: We have to analyze $L^{n-2} \to F'$. This morphism is
injective
and completes to the exact sequence:

$$
0 \to \frac{J_{n-2}}{J_{n-1}} \to \frac{I_{n-2}}{I_{n-1}} \to
\frac{I_{n-2}}{J_{n-2}+I_{n-1}} \to 0 \ .
$$
As the surjective morphism between locally free rank $1$ sheaves,
$I_{n-2}/(J_{n-2}+I_{n-1}) \to I_2/J_2$
must be an isomorphism, the  above exact sequence is in fact:
$$
0 \to L^{n-2} \to F' \to K \to 0
$$

\noindent Dualizing this sequence and tensoring with $L^n$, one obtains:

$$
0 \to K^\vee \otimes L^n \to F \to L^2 \to 0
$$

\noindent Comparing with an exact sequence from above, one gets $K^\vee \otimes L^n \cong
K^2$, i.e.
$$
K^3 \cong L^n \ .
$$

\noindent {\em Degree $n-1$}: We have to analyze the morphism $L^{n-1} \to E' $. This
completes to the exact sequence:
$$
0 \to J_{n-1}/J_n 
\to I_{n-1}/I_n \to I_{n-1}/(J_{n-1}+I_n) \to 0 \ ,
$$
i.e. :
$$
0 \to L^{n-1} \to E' \to K^2 \to 0
$$

\noindent which, dualizing and tensoring with $L^n$ gives the exact sequence:
$$
0 \to K^{-2}\otimes L^n \to E \to L \to 0 \ .
$$
\noindent One obtains again $K^3 \cong L^n$.

\medskip
\noindent Conversely:
\begin{theorem}
Let $X \subset P$ be two smooth varieties, $\hbox{codim}_PX \ge2$. Let $I$ be the
sheaf-ideal of $X$ in $P$. All cuspidal nilpotent structures of type $C_{3,n}$ on $X$,
embedded in $P$ can be obtained in the following way:\\
{\bf Step 1}. Give two line bundles $L$, $K$ on $X$, satisfying $L^n \cong K^3$ and
an extension
$0 \to K \stackrel{\iota}{\longrightarrow} E \stackrel{\pi}{\longrightarrow}L \to 0$ .
The exact sequence gives rise to an injection $K^2 \hookrightarrow S^2E$ and a
surjection $S^2E \to L^2$. Give an extension 
  $0 \to K^2 \to F \to L^2 \to 0$ and a  surjection $S^2E \to F$ such that the second
extension fits in the commutative diagram:
$$
\xymatrix{
& & S^2E \ar[d] \ar[rd] & & \\
0 \ar[r] & K^2 \ar[r] \ar@{^{(}->}[ru] &  F  \ar[r]  &
L^2 \ar[r]  & 0 \\
}
$$
{\bf Step 2}. Give two surjections $p_2:I/I^2 \to E$, $q_2:I/I^2 \to L$, such that the
following diagram is commutative:

$$
\xymatrix{
I/I^2  \ar[rr]^{p_2}\ar@{=}[d]& & E \ar[r] \ar[d] & 0 \\
I/I^2  \ar[rr]^{q_2} & & L \ar[r] & 0 \\
}
$$
\medskip

\noindent Take $J_2 = \hbox{ker}(I \to I/I^2 \to E)$, $I_2
=\hbox{ker}(I \to I/I^2 \to L)$. Then $K \cong I_2/J_2$, $K\otimes E \cong II_2/IJ_2$,
and, in convenient local coordinates $J_2 = (x^2, xy, y^2, \bz)$, $I_2 = (x^2, y, \bz)$.

\noindent {\bf Step 3.} Observe that: $ S^2E \cong I^2/IJ_2 \hookrightarrow J_2/IJ_2$.
Give surjections $p_3: J_2/IJ_2 \to F$, $q_3:I_2/II_2 \to L^2$ such that the following
diagram is commutative:

$$
\xymatrix{
I_2^2/I_2J_2 = K^2 \ar@{_{(}->}[dd] \ \ \ \ \ar@{^{(}->}[rd]& &\\
& S^2E \cong I^2/IJ_2 \ar@{_{(}->}[ld] \ar[rd]& & \\
J_2/IJ_2 \ar[rr]^{p_3}\ar[d]& & F \ar[r] \ar[d] & 0 \\
I_2/II_2 \ar[rr]^{q_3} & & L^2 \ar[r] & 0 \\
}
$$
\medskip
Take $J_3 =\hbox{ker}( J_2 \to J_2/IJ_2 \to F)$, $I_3 =\hbox{ker}(I_2 \to I_2/II_2 \to
L^2)$. In convenient local coordinates: $J_3=(x^3,xy,y^3, \bz)$, $I_3=(x^3,y,\bz)$

\noindent {\bf Step 4}. Take $q_4$ a retract of the canonical injection $L^3 \cong
I^3/I^2I_2 \hookrightarrow I_3/II_3$ and a surjection $p_4: J_3/IJ_3 \to L^3 $ such that
the following  digram is commutative:

$$
\xymatrix{
J_3/IJ_3 \ar[rr]^{p_4}\ar[d]& & L^3 \ar[r] \ar@{=}[d] & 0 \\
I_2/II_2 \ar[rr]^{q_4} & & L^3 \ar[r] & 0 \\
}
$$
\medskip
Take $J_4 =\hbox{ker}( J_3 \to J_3/IJ_3 \to L^3)$, $I_4 =\hbox{ker}(I_3 \to I_3/II_3 \to
L^3)$. Then, in convenient local coordinates: $J_4 =(x^4,xy,y^3,\bz)$, $I_4 =(x^4,y,\bz)$.

\medskip
\hfil \vdots
\medskip

\noindent {\bf Step $\ell$, $\ell \le n-2$}. With $J_{\ell -1}=(x^{\ell -1}, xy,y^3,\bz)$,
$I_{\ell -1}=(x^{\ell -1},y,\bz)$, take $q_\ell$ a retract of the canonical injection
$L^{\ell -1} \cong I^{\ell -1}/I^{\ell -2}I_2 \hookrightarrow I_{\ell -1}/II_{\ell -1}$
and $p_\ell$ such that the following diagram is commutative:
$$
\xymatrix{
J_{\ell -1}/IJ_{\ell -1} \ar[rr]^{p_\ell }\ar[d]& & L^{\ell -1} \ar[r] \ar@{=}[d] & 0 \\
I_{\ell -1}/II_{\ell -1} \ar[rr]^{q_\ell } & & L^{\ell -1} \ar[r] & 0 \\
}
$$
\medskip
It follows: $J_\ell =(x^\ell , xy, y^3,\bz)$, $I_\ell =(x^\ell ,y,\bz)$.

\medskip
\hfil \vdots
\medskip

\noindent {\bf Step n-1}. From the previous step we have $J_{n-2}=(x^{n-2},xy,y^3,\bz)$,
$I_{n-2}=(x^{n-2},y,\bz)$. Dualizing the extension which gives $F$ and tensoring with
$L^n$, one gets the exact sequence:
$$
0 \to L^{n-2}\to F' \to K \to 0 \ .
$$
Take $p_{n-1}: J_{n-2}/IJ_{n-2} \to L^{n-2}$ to be a retract of the canonical embedding
$L^{n-2}\cong I^{n-2}/I^{n-3}I_2 \hookrightarrow J_{n-2}/IJ_{n-2}$ and $q_{n-2}$ such that
the following diagram is commutative:
$$
\xymatrix{
J_{n-2}/IJ_{n-2} \ar[rr]^{p_{n-1} }\ar[d]& & L^{n-2} \ar[r] \ar@{^{(}->}[d] & 0
\\
I_{n-2}/II_{n-2} \ar[rr]^{q_{n-1} } & & F' \ar[r] & 0 \\
}
$$
\medskip
Take $J_{n-1}=\hbox{ker}(J_{n-2} \to J_{n-2}/IJ_{n-2} \to L^{n-2})$ and $I_{n-1}= \\
\hbox{ker}(I_{n-2} \to I_{n-2}/II_{n-2} \to L^{n-2})$.
It follows that, in convenient local coordinates: $J_{n-1}=(x^{n-1},xy,y^3,\bz)$ and
$I_{n-1}=(x^{n-1},xy,y^2.x\bz ,y\bz ,\bz^2)$\medskip

\noindent One shows easily:
$$
\frac{I_{n-2}}{J_{n-2}+II_{n-2}} \cong K
$$

\noindent {\bf Step n}. Take $p_n: J_{n-1}/IJ_{n-1} \to L^{n-1}$ to be a retract of
the canonical embedding $L^{n-1}\cong I^{n-1}/I^{n-1}I_2 \hookrightarrow
J_{n-1}/IJ_{n-1}$. Take $J_n = I_n =
\hbox{ker}(J_{n-1}\to J_{n-1}/IJ_{n-1} \to L^{n-1})$. Then, in convenient local
coordinates : $J_n=I_n \cong (x^n,xy,y^3, x\bz , y\bz , \bz ^2)$ and  $E' \cong
I_{n-1}/I_n$.\medskip

\noindent {\bf Step n+1}. Take $p_n =q_n : I_n/II_n \to L^3 \cong K^3 $ to be a retract of
the canonical injections: $L^n \cong I^n/I^{n-1}I_2 \to J_n/IJ_n $ and $K^3 \cong
I_2^2I_{n-2}/(I_2J_2I_{n-2}+I_2^2J_{n-2}+II_2^2I_{n-1}) \to I_n/II_n$. Then
$J_{n+1}=I_{n+1} =
\hbox{ker}(J_n \to J_n/IJ_n \to L^n)$ is a nilpotent structure on X of type $C_{3,n}$.

\medskip
\end{theorem}

\Proof Like in the previous case, one has to verify step by step all assertions.

\qed

\bigskip

{\bf Acknowledgement}. In the preparation  of this paper the author had no
other support besides  the membership to the Institute of Mathematics of the Romanian
Academy.

\noindent Nicolae Manolache\\
Institute of Mathematics "Simion Stoilow"\\
of the Romanian Academy \\
P.O.Box 1-764
Bucharest, RO-014700

\noindent e-mail: nicolae.manolache@imar.ro
\end{document}